\def \QQ {\mathbb Q}

\def \epsilon{\varepsilon}

\def \OO  {{\mathcal O}}

\def \d {\text{d}}

\def \fine {{\hfill \qedsymbol}}

\documentclass[12pt,reqno]{amsart}

\usepackage{amsmath,amssymb}

\usepackage{url} 

\begin{document}


\title[]{Primes and prime ideals in short intervals}
\author[]{L.GRENI\'E, G.MOLTENI \lowercase{and} A.PERELLI}
\maketitle

{\bf Abstract.} We prove the analog of Cram\'er's short intervals theorem for primes in arithmetic progressions and prime ideals, under the relevant Riemann Hypothesis. Both results are uniform in the data of the underlying structure. Our approach is based mainly on the inertia property of the counting functions of primes and prime ideals.



\section{Introduction}

\smallskip
A famous theorem of H.Cram\'er \cite{Cra/1920} states that, assuming the Riemann Hypothesis, there is always a prime between $x$ and $x+h$ provided $x\geq x_0$ and $c_1\sqrt{x}\log x\leq h \leq x$, with suitable constants $x_0,c_1>0$. Actually, under the same assumptions we have that
\[
\pi(x+h) - \pi(x) \geq c_2 \frac{h}{\log x}
\]
with a suitable $c_2>0$, and also that 
\[
\pi(x+h) - \pi(x) \sim \frac{h}{\log x}
\]
provided $x\geq h = \infty(\sqrt{x}\log x)$. Here $f(x)=\infty (g(x))$ means that $f(x)/g(x) \to \infty$ as $x\to\infty$. Apart from the explicit value of the involved constants, this is still the best known result about primes in short intervals, under the Riemann Hypothesis. Sharper results can be obtained assuming various forms of the pair-correlation conjecture for the zeta zeros, see {\it e.g.} Heath-Brown \cite{HB/1982}, Languasco-Perelli-Zaccagnini \cite{LPZ/I} and the literature quoted there. A simple proof of Cram\'er's theorem can be obtained from a suitable smoothed explicit formula for $\psi(x)$, see the footnote on p.256 of Ingham \cite{Ing/1937}.

\smallskip
In this paper we show that rather general theorems of Cram\'er's type follow, under the appropriate Riemann Hypothesis, from two results often ready in the literature, namely a short intervals mean-square estimate and a Brun-Titchmarsh type theorem. Indeed, the latter result implies that the relevant counting function satisfies a suitable inertia property, which is then played against the short intervals mean-square bound to get a contradiction if the interval is not too short. We illustrate our approach in the case of primes in arithmetic progressions and of prime ideals, since apparently these results do not appear in the literature. In the first case all the ingredients are already known, so we proceed directly to the proof of Cram\'er's theorem for arithmetic progressions. In the case of algebraic number fields we first deal with the required ingredients; see in particular Proposition 1 below, which is of some independent interest. In both cases our results are uniform in the data of the underlying structure. However, in the second case the inertia method gives Proposition 3, which in the uniformity aspect is weaker than Theorem 2, proved here by the classical smoothed explicit formula approach. This is due to the lack, in the current literature, of sharp uniform bounds of Brun-Titchmarsh type for number fields. We shall discuss this issue later on in the paper.

\smallskip
As usual, for $(a,q)=1$ we write
\[
\pi(x;q,a) = \sum_{\substack{p\leq x \\ \text{$p\equiv a$ (mod $q$)}}} 1
\]
and let $\varphi(q)$ denote Euler's function. Moreover, given an algebraic number field $K$ of degree $n_K$, we denote 
by $d_K$ the absolute value of its discriminant, by $\frak{P}$ the prime ideals of the ring $\OO_K$ of the integers of $K$, by $N(\frak{P})$ their norm and write
\[
\pi_K(x) = \sum_{N(\frak{P})\leq x} 1.
\]
Finally, given an integer $q\geq 1$ and a number field $K$, we denote by GRH and DRH the Riemann Hypothesis for the Dirichlet $L$-functions associated with the characters $\chi$ (mod $q$) and for the Dedekind zeta function $\zeta_K(s)$, respectively. With this notation, our main results are as follows.

\medskip
{\bf Theorem 1.} {\sl Let $(a,q)=1$ and assume GRH. Then there exist absolute constants $x_0,c_1,c_2>0$ such that for $x\geq x_0$ and $c_1 \varphi(q)\sqrt{x}\log x \leq h \leq x$ we have}
\[
\pi(x+h;q,a)-\pi(x;q,a) \geq c_2 \frac{h}{\varphi(q)\log x}.
\]

\medskip
Clearly, under the same assumptions the same argument also gives
\[
\pi(x+h;q,a)-\pi(x;q,a) \sim \frac{h}{\varphi(q)\log x}
\]
provided $x\geq h =\infty(\varphi(q)\sqrt{x}\log x)$.

\medskip
{\bf Theorem 2.} {\sl Assume DRH for the number field $K$. Then there exist absolute constants $x_0,c_1,c_2>0$ such that for $x\geq x_0$ and $c_1 (n_K\log x +\log d_K)\sqrt{x} \leq h \leq x$ we have}
\[
\pi_K(x+h)-\pi_K(x) \geq c_2 \frac{h}{\log x}.
\]

\medskip
As before, the same proof shows also that
\begin{equation}
\label{1}
\pi_K(x+h)-\pi_K(x) \sim \frac{h}{\log x}
\end{equation}
provided $x\geq h =\infty\big((n_K\log x +\log d_K)\sqrt{x}\big)$. Note that Theorem 2 represents an instance of Lang's \cite{Lan/1971} ``recipe'' asserting that, broadly speaking, when extending to a number field $K$ the classical results known for $\QQ$ one should replace $\log x$ by $n_K\log x+\log d_K$. Note also that if $K$ is a cyclotomic field then the quality of the $K$-uniformity in Theorem 2 is comparable to the $q$-uniformity in Theorem 1.

\smallskip
We conclude remarking that the technique in the proof of the above theorems works for rather general counting functions, giving individual short intervals results as soon as suitably sharp short intervals mean-values and inertia type result are available.

\section{Proofs}

\smallskip
As customary, we prove Theorem 1 for the $\psi$-function and then the required result is recovered by elementary arguments since $h/\varphi(q)$ is large enough. Let $X$ be sufficiently large, $q,h\leq X$, $(a,q)=1$ and write 
\[
\Delta(x,h) = \psi(x+h;q,a)-\psi(x;q,a)-\frac{h}{\varphi(q)}.
\]
The required mean-square bound follows from a result of Prachar \cite{Pra/1976} under GRH (see also 
\eject
\noindent
Goldston-Y\i ld\i r\i m \cite{Go-Yi/1998}), namely
\begin{equation}
\label{2}
\int_X^{2X} |\Delta(x,h)|^2dx \ll hX\log^2(qX),
\end{equation}
where the constant in the $\ll$-symbol is absolute. Let now $h/\varphi(q)>X^{1/10}$. From the well known Brun-Titchmarsh theorem, see Montgomery-Vaughan \cite{Mo-Va/1973}, we deduce that if there exists $\overline{x}\in (X,2X)$ such that
\begin{equation}
\label{3}
|\Delta(\overline{x},h)| > \frac{1}{4} \frac{h}{\varphi(q)},
\end{equation}
then
\begin{equation}
\label{4}
|\Delta(x,h)| > c \frac{h}{\varphi(q)}
\end{equation}
for all $x\in(\overline{x}-c'h,\overline{x}+c'h)$, with certain absolute constants $c,c'>0$. Inequalities \eqref{3} and \eqref{4} express the  inertia property of the $\psi$-function (see also Theorem 1 of Bazzanella-Perelli \cite{Ba-Pe/2000}).

\smallskip
Let now
\[
E(X,h) = \{x\in[X,2X]: |\Delta(x,h)| >\frac{1}{4} \frac{h}{\varphi(q)}\}
\]
and suppose that $E(X,h) \neq \emptyset$. Then from \eqref{2}-\eqref{4} we get
\[
h (\frac{h}{\varphi(q)})^2 \ll \int_X^{2X} |\Delta(x,h)|^2dx \ll hX\log^2(qX),
\]
thus $h \ll \sqrt{X} \varphi(q)\log(qX)$. Hence, with suitable absolute constants in the $\gg$-symbols and provided $x$ is sufficiently large,
\[
\psi(x+h;q,a)-\psi(x;q,a) \gg  \frac{h}{\varphi(q)}
\]
if $x \geq h\gg \varphi(q)\sqrt{x}\log(qx)$. Theorem 1, and the statement after it, therefore follow. \fine

\medskip
As anticipated in the Introduction, in the number fields case we first present the proof of a weaker form of Theorem 2 in the uniformity aspect, obtained by the inertia approach. Write
\[
\psi_K(x) = \sum_{N(\frak{P}^m)\leq x} \log N(\frak{P})=\sum_{n\leq x}\Lambda_K(n) \ \ \text{say, and} \ \ \Delta_K(x,h) = \psi_K(x+h)-\psi_K(x)-h,
\]
and let $L=\log X$. The analog of \eqref{2} is given by

\medskip
{\bf Proposition 1.} {\sl Assume DRH for the number field $K$. Then there exist absolute constants $c,X_0>0$ such that for $X\geq X_0$ and $2\leq h\leq X$ we have}
\[
\int_{X}^{2X} |\Delta_K(x,h)|^2\d x \leq c X(h+L^2)(n_KL+\log d_K)^2.
\]

\medskip
{\it Proof.} Let $2\leq T \leq x$; the constants in the $O$- and $\ll$-symbols below are absolute. Denoting by $N_K(T)$ the number of zeros $\rho=\beta+i\gamma$ of $\zeta_K(s)$ with $0\leq \beta\leq 1$ and $|\gamma|\leq T$, using the notation in the Introduction we unconditionally have
\begin{equation}
\label{5}
N_K(T) = \frac{n_K}{\pi} T\log T + \frac{T}{\pi} \log\left(\frac{d_K}{(2\pi e)^{n_K}}\right) + O\big(\log(d_KT^{n_K})\big);
\end{equation}
see Kadiri-Ng \cite{Ka-Ng/2012}. Moreover, using \eqref{5} in the explicit formula in Lagarias-Odlyzko \cite{La-Od/1977} (specialized to the case of $\zeta_K(s)$) we have, again unconditionally, that
\begin{equation}
\label{6}
\psi_K(x) = x - \sum_{|\gamma|\leq T} \frac{x^\rho}{\rho} +R_K(x,T)
\end{equation}
with
\[
R_K(x,T)  =O\left(\frac{x}{T}(n_K\log x + \log d_K)\log x \right).
\]
Proposition 1 follows now from \eqref{4} and \eqref{6} by the classical arguments in Lemmas 5 and 6 of Saffari-Vaughan \cite{Sa-Va/1977} (notice a misprint in (6.20) there, where $h^2$ on the r.h.s. should be replaced by $h$); here is a brief sketch. Arguing as in Lemma 6 of \cite{Sa-Va/1977} we have (see (6.21) there)
\begin{equation}
\label{7}
\int_{X}^{2X} |\Delta_K(x,h)|^2\d x  \ll \frac{X}{h} \int_{h/3X}^{3h/X} \left(\int_X^{3X} |\psi_K(x+\theta x)-\psi_K(x)-\theta x|^2 \d x\right) \d \theta.
\end{equation}
Choosing $T=X$ in \eqref{6}, the contribution of $R_K(x,T)$ to the r.h.s. of \eqref{7} is 
\begin{equation}
\label{8}
\ll XL^2(n_KL+\log d_K)^2,
\end{equation}
while the contribution of the remaining part of the explicit formula is, thanks to \eqref{5},
\begin{equation}
\label{9}
\ll X^2 \big(\frac{h}{X}\big)^2 N_K\left(\frac{X}{h}\right) \max_{2\leq t\leq X/h} \big(N_K(t+1)-N_K(t)\big) \ll hX(n_KL+\log d_K)^2.
\end{equation}
Proposition 1 follows from \eqref{7}-\eqref{9}. \fine

\medskip
Proposition 1 represents another instance of Lang's ``recipe'' reported in the Introduction. As far as we know, such a phenomenon has not been established in the case of Brun-Titchmarsh type bounds, and actually it is not clear to us how should the right extension look like in this case; we briefly discuss this issue at the end of the section. Hence we use the following simple but uniform bound, which however is unlikely to be sharp in the range needed here. 

\medskip
{\bf Proposition 2.} {\sl Let $K$ be a number field and $2\leq h\leq x$. Then}
\begin{equation}
\label{10}
\pi_K(x+h) - \pi_K(x) \leq 4n_K \frac{h}{\log h}.
\end{equation}

\medskip
{\it Proof.} Again we use the notation in the Introduction. Let $\{k_j\}_{j\in J}$, $k_j\geq 1$, be the exponents of the prime powers in the interval $(x,x+h]$; clearly, $|J|\leq h+1$. Since it is well known that for $1\leq k\leq n_K$ there are at most $n_K/k$ prime ideals of $\OO_K$ with norm $p^k$, we have
\[
\pi_K(x+h)-\pi_K(x) \leq n_K \sum_{j\in J}\frac{1}{k_j} \big( \pi((x+h)^{1/k_j}) - \pi(x^{1/k_j})\big).
\]
But $(x+h)^{1/k} - x^{1/k} \leq x^{1/k}h/(kx)$, hence, applying to $\pi((x+h)^{1/k_j}) - \pi(x^{1/k_j})$ the Brun-Titchmarsh theorem when $k_j=1$ (Montgomery-Vaughan \cite{Mo-Va/1973} with modulus $q=1$) and the trivial bound $\leq h/(k_j\sqrt{x})+1$ when $k_j\geq 2$, we get
\[
\pi_K(x+h)-\pi_K(x) \leq 2n_K \frac{h}{\log h} + n_K \sum_{j\in J, k_j\geq2} \frac{1}{k_j}\big(\frac{h}{\sqrt{x}} \frac{1}{k_j} +1\big).
\]
Since clearly $\sum_{j\in J, k_j\geq 2} 1/k_j \leq \log(|J|+1) \leq \log(h+2)$, Proposition 2 follows by a simple computation. \fine

\medskip
{\bf Proposition 3.} {\sl Theorem $2$ holds with $c_1 n_K^{1/2}(n_K\log x +\log d_K)\sqrt{x} \leq h \leq x$ in place of $c_1 (n_K\log x +\log d_K)\sqrt{x} \leq h \leq x$.}

\medskip
{\it Proof.} We argue along the lines of Theorem 1. Indeed, for $X$ sufficiently large and {e.g.} $h/n_K > X^{1/10}$, from Proposition 2 we have that if there exists $\overline{x}\in (X,2X)$ with
\[
|\Delta_K(\overline{x},h)| > \frac{1}{4} h, \quad \text{then} \quad |\Delta_K(x,h)| > c h
\]
for all $x\in(\overline{x}-c'h/n_K,\overline{x}+c'h/n_K)$, with certain absolute constants $c,c'>0$. Playing this against Proposition 1 we therefore obtain that
\eject
\[
\frac{h^3}{n_K} \ll (n_K L+\log d_K)^2 hX.
\]
Hence $h \ll n_K^{1/2}(n_KL + \log d_K)\sqrt{X}$, and Proposition 3 follows. \fine

\medskip
The sharper result stated in Theorem 2 is obtained using the direct approach by the smoothed explicit formula. We follow the general lines of the proofs in Dudek \cite{Dud/2015} and Dudek-Greni\'e-Molteni \cite{D-G-M/2015}, where explicit versions of Ingham's approach to Cram\'er's theorem are developed. Integrating the infinite explicit formula for $\psi_K(x)$ from 2 to $x$, see (1.3a) and Lemmas 3.2 and 3.3 in Greni\'e-Molteni \cite{Gr-Mo/2015}, we obtain
\[
\int_2^x\psi_K(t)\d t = \frac{x^2}{2} - \sum_\rho \frac{x^{\rho+1}}{\rho(\rho+1)} -c_Kx + c'_K  + O(n_K x \log x),
\]
where $\rho$ runs over the non-trivial zeros of $\zeta_K(s)$ and $c_K,c'_K$ are certain constants depending on $K$; we are not concerned with their value since $c_K$ and $c'_K$ simply disappear after the manipulations leading to the next displayed equation. Introducing the weights $w(n)=\max(1-\frac{|x-n|}{h},0)$ as in the last row of p.773 of \cite{Dud/2015} and arguing as on p.774 there we get
\[
\sum_{x-h<n<x+h} \Lambda_K(n)w(n) = h - \frac{1}{h} \sum_{\rho}\frac{(x+h)^{\rho+1}-2x^{\rho+1}+(x-h)^{\rho+1}}{\rho(\rho+1)} +  O\Big(n_K \frac{x}{h}\log x\Big).
\]
Now we split the sum over the $\zeta_K$-zeros into the subsums $\Sigma_1$ and $\Sigma_2$ cutting at $T=x/h$, and use DRH and \eqref{5} as in the proof of Theorem 1.2 of \cite{Dud/2015}, thus obtaining
\begin{equation}
\label{11}
W(x,h) := \sum_{x-h<n<x+h} \Lambda_K(n)w(n) =  h+ O((n_K\log x + \log d_K)\sqrt{x}) + O\Big(n_K \frac{x}{h}\log x\Big).
\end{equation}
From \eqref{11} we obtain the behavior of the unweighted sum, observing that for every $0<\epsilon<1$
\[
\begin{split}
-\frac{1}{\epsilon}\big((1-\epsilon)W(x,(1-\epsilon)h) &- W(x,h)\big) \leq \psi_K(x+h)-\psi_K(x-h) \\
&\leq \frac{1}{\epsilon}\big((1+\epsilon)W(x,(1+\epsilon)h) - W(x,h)\big)
\end{split}
\]
since $\Lambda_K(n)\geq0$. Theorem 2 and the assertion after it follow at once. \fine

\medskip
We conclude with a brief discussion on the Brun-Titchmarsh theorem for number fields and its relevance to this paper. Note that the dependence on the data of $K$ in Proposition 2, where bounded $h$ are allowed, is essentially best possible. Indeed, if a prime $p\in(x,x+2]$, say, splits in $O_K$ into the product of $n_K$ prime ideals of norm $p$, then clearly $\pi_K(x+2)-\pi_K(x) \geq n_K$, while Proposition 2 gives $\pi_K(x+2)-\pi_K(x) \leq cn_K$ with some absolute $c>0$. Note that, although the constant in the classical Brun-Titchmarsh theorem is of great interest, the absolute constant in front of $n_K$ in \eqref{10} plays essentially no role in this paper.  For larger $h$ the dependence on $K$ in \eqref{10} is unsatisfactory, as indeed the prime ideal theorem, or \eqref{1}, shows. 

\smallskip
The bounds of Brun-Titchmarsh type are usually obtained by the Selberg sieve. Apparently, an application of the Selberg sieve to $\pi_K(x+h)-\pi_K(x)$, see e.g. Hinz-Loedemann \cite{Hi-Lo/1994}, brings into play the residue $\nu_K$ of the Dedekind zeta function $\zeta_K(s)$. It is well known that $\nu_K$ depends on several invariants of $K$, and even under DRH its dependence on such invariants is not completely under control. This adds some difficulties to the problem of obtaining sharp versions of Proposition 2. Perhaps one can prove that
\[
\pi_K(x+h)-\pi_K(x) \leq c h/\log(h/d_K),
\]
but this is weaker than what is obtainable for an abelian extension $K/\QQ$, namely
\eject
\begin{equation}
\label{12}
\pi_K(x+h)-\pi_K(x) \leq c h/\log(h/q_K),
\end{equation}
where $q_K$ is the conductor of $K$. Bound \eqref{12} can be obtained coupling the classical Brun-Titchmarsh theorem for arithmetic progressions with the Kronecker-Weber theorem for abelian extensions of $\QQ$. Actually, when \eqref{12} is coupled with Proposition 1 we get back, in the abelian case, a result of the same quality as Theorem 1.

\medskip
{\bf Acknowledgements.} We wish to thank Olivier Ramar\'e for detecting some inaccuracies in a previous version. The authors are members of the INdAM groups GNSAGA and GNAMPA.


\ifx\undefined\bysame{poly}.
\newcommand{\bysame}{\leavevmode\hbox to3em{\hrulefill}\ ,}
\fi

\bigskip
\noindent
Lo\"ic Greni\'e, Dipartimento di Ingegneria Gestionale, dell'Informazione e della Produzione, Universit\`a di Bergamo, viale Marconi 5, 24044 Dalmine (BG), Italy. e-mail: loic.grenie@gmail.com

\smallskip
\noindent
Giuseppe Molteni, Dipartimento di Matematica, Universit\`a di Milano, via Saldini 50, 20133 Milano, Italy. e-mail:
giuseppe.molteni1@unimi.it

\smallskip
\noindent
Alberto Perelli, Dipartimento di Matematica, Universit\`a di Genova, via Dodecaneso 35, 16146 Genova, Italy. e-mail: perelli@dima.unige.it

\end{document}